\documentclass[11pt]{amsart}

\usepackage{amssymb}
\usepackage{amsthm}
\usepackage{mathrsfs}
\usepackage{cite}
\usepackage{upref}
\usepackage{verbatim}
\usepackage{xspace}

\renewcommand{\bar}{\overline}
\newtheorem{theorem}{Theorem}

\newtheorem{conjecture}[theorem]{Conjecture}

\theoremstyle{definition}

\newtheorem{remark}[theorem]{Remark}

\numberwithin{theorem}{section}
\numberwithin{equation}{section}

\def\R{\mathbb R}
\def\C{\mathbb C}
\def\Q{\mathbb Q}
\def\Z{\mathbb Z}
\def\G{\mathbf{G}}
\def\fg{\mathfrak{g}}
\def\fk{\mathfrak{k}}
\def\E{\mathscr{E}}

\def\Eis{\text{Eis}}
\def\cusp{\text{cusp}}
\def\grit{\text{G}}
\def\nongrit{\text{nG}}
\def\disc{\text{disc}}
\DeclareMathOperator\Frob{{Frob}}
\DeclareMathOperator\Gal{{Gal}}
\def\Hom{\text{Hom}}
\newcommand{\GL}{{\rm GL}}

\newcommand{\Sp}{{\rm Sp}}
\newcommand{\SO}{{\rm SO}}

\newcommand{\SL}{{\rm SL}}

\newcommand{\thmcite}[1]{\emph{\cite{#1}}}

\newcommand{\gbg}{\Gamma \backslash  G}
\newcommand{\Vor}{Vorono\v{\i}\xspace}

\title[On the cohomology of congruence subgroups of $\SL_{4} (\Z)$]{On the cohomology of congruence subgroups of $\SL_{4} (\Z )$}
\author{Paul~E. Gunnells}
\address{Department of Mathematics and Statistics,
University of Massachusetts,
Amherst, MA 01003, USA
}
\email{gunnells@math.umass.edu}

\date{May 4, 2009}
\thanks{Partially supported by the NSF}

\begin{document}

\begin{abstract}
We survey our joint work with Avner Ash and Mark McConnell that
computationally investigates the cohomology of conguence subgroups of
$\SL_{4} (\Z)$.
\end{abstract}

\maketitle

\section{Introduction}

Let $\G$ be a semisimple algebraic group defined over $\Q$ and let
$\Gamma \subset \G (\Q)$ be an arithmetic group.  The cohomology $
H^{*} (\Gamma ; \C) $ plays an important role in number theory, in
that it provides a concrete realization of certain automorphic forms
(\S\ref{s:coho-autfms}).  For instance, if $\G = \SL_{2}$ and $\Gamma$
is a congruence subgroup, then by the Eichler--Shimura isomorphism the
cohomology space $H^{1} (\Gamma ; \C)$ is built from weight two
modular forms of level $N$ (\S\ref{ss:modularforms}).  Not every
automorphic form can appear in the cohomology of an arithmetic group,
but those that do are widely believed to be related to arithmetic
geometry, such as Galois representations and motives
(\S\ref{s:arithmetic}).  The cohomology and its structure as a Hecke
module can be computed very explicitly in many cases.  Thus the
cohomology of arithmetic groups provides a tool to formulate and test
conjectures about the links between automorphic forms and arithmetic.

In a series of paper with Avner Ash and Mark McConnell \cite{agm1,
agm2, agm3}, we have investigated the group cohomology for $\G =
\SL_{4}$ and $\Gamma = \Gamma_{0} (N)$ the congruence subgroup of
$\SL_{4} (\Z)$ with bottom row congruent to $(0,0,0,*) \bmod N$.  Our
work has focused on a particular cohomology space attached to
$\Gamma$, namely $ H^{5} (\Gamma ; \C)$.  (The choice of degree $5$ is
explained below.)  We have computed the dimension of this cohomology
space for prime levels $N\leq 211$, and in many cases have computed
the action of the Hecke operators.

The purpose of this note is to give an introduction to \cite{agm1,
agm2, agm3} and to explain how we perform these computations
(\S\ref{s:techniques}).  We also discuss our computational results
(\S\ref{s:results}), which indicate relationships between $ H^{5}
(\Gamma ; \C)$ and elliptic modular forms, the cohomology of subgroups
of $\SL_{3} (\Z)$, and certain Siegel modular forms, the paramodular
forms.

\subsection*{Acknowledgements.} This article is based on a lecture
delivered by the author at the 2009 RIMS conference \emph{Automorphic
representations, automorphic $L$-functions and arithmetic}.  The
author thanks the organizers of the conference, Yoshi-Hiro Ishikawa
and Masao Tsuzuki, for the opportunity to speak.  He also warmly
thanks Tomoyoshi Ibukiyama for the invitation to Japan and for being
such an excellent host.  Finally, he thanks Avner Ash and Mark
McConnell for many years of interesting and fun mathematics.

\section{Cohomology and automorphic forms}\label{s:coho-autfms}
In this section we briefly explain what the cohomology of arithmetic
groups has to do with automorphic forms.  For more information we
refer to \cite{borel.wallach, borel.survey, schwermer.survey,
li.schwermer.survey, vogan, stein.book, harder.book}.

\subsection{}\label{ss:notation} Let $\G$ be a semisimple connected
algebraic group defined over $\Q$, and let $\Gamma \subset \G (\Q)$ be
an arithmetic subgroup.  Let $E$ be a finite-dimensional rational
complex representation of $\G(\Q)$.  Then $E$ is naturally a $\Z
\Gamma$-module, and we are interested in the group cohomology $H^{*}
(\Gamma ; E)$.

We can compute $H^{*} (\Gamma ; E)$ topologically as follows.  Let $G
= \G (\R)$ be the group of real points of $\G$; $G$ is a semisimple
Lie group.  Let $K \subset G$ be a maximal compact subgroup, and let
$X = G/K$ be the associated global symmetric space.  The space $X$ is
contractible with a left $\Gamma$-action.  If $\Gamma$ is
torsion-free, then $\Gamma \backslash X$ is an Eilenberg--Mac~Lane
space, and we have
\begin{equation}\label{eqn:gpcohoiso}
H^{*} (\Gamma ; E) = H^{*} (\Gamma \backslash X; \E), 
\end{equation}
where $\E$ is the local coefficient system on the manifold $\Gamma
\backslash X$ attached to $E$.  In fact, since we consider only
complex representations $E$, the isomorphism \eqref{eqn:gpcohoiso}
holds even if $\Gamma$ has torsion.  In this case the quotient $\Gamma
\backslash X$ is an orbifold, but nevertheless one can construct a
local system $\E$ on the quotient such that \eqref{eqn:gpcohoiso}
remains true.

\subsection{} The first step in getting automorphic forms into the
picture is the de Rham theorem.  Let $\Omega^{p} = \Omega^{p} (X, E)$
be the space of $E$-valued $p$-forms on $X$, and let $\Omega^{p} (X,
E)^{\Gamma}$ be the subspace of $\Gamma$-invariant forms. We have a
differential $d\colon \Omega^{p}\rightarrow \Omega^{p+1}$ and an
isomorphism of cohomology spaces \cite[Theorem 2.2.2]{borel.wallach}
\[
H^{*} (\Gamma ; E) = H^{*} (\Omega^{*} (X,E)^{\Gamma}).
\]

Now recall that $X=G/K$.  Let $\fg$ (respectively $\fk$) be the Lie
algebra of $G$ (resp. $K$).  We can identify the tangent space to $X$
at the basepoint with the quotient $\fg / \fk$.  Since differential
forms are sections of the exterior powers of the tangent bundle, one
can show \cite[Prop.~1.5]{vogan}
\[
\Omega^{p} (\Gamma \backslash X, \C) = \Hom_{K} (\wedge^{p}(\fg
/\fk), C^{\infty} (\Gamma \backslash G)),
\]
where on the right we take complex-valued smooth functions on $\Gamma
\backslash G$, and where $K$ acts on $\wedge^{p}(\fg /\fk)$ by the
adjoint action and on $C^{\infty} (\Gamma \backslash G)$ by right
translations.  More generally, if we want cohomology with
coefficients, then we take $E$-valued differential forms, and we have
\begin{equation}\label{eqn:derham}
\Omega^{p} (\Gamma \backslash X, E) = \Hom_{K} (\wedge^{p}(\fg
/\fk), C^{\infty} (\Gamma \backslash G)\otimes E).
\end{equation}
The right hand side of \eqref{eqn:derham} inherits a differential from the left.
The resulting cohomology spaces are denoted
\[
H^{*} (\fg , K; C^{\infty} (\Gamma \backslash G)\otimes E).
\]
This cohomology is called the \emph{relative Lie algebra cohomology of
$(\fg , K)$ with coefficients in $E$}, or simply \emph{$(\fg ,
K)$-cohomology}.

\subsection{}
To summarize, we have identifications
\begin{equation}\label{eqn:gKeqn}
H^*(\Gamma ; E) = H^{*} (\Gamma \backslash X; \E )= H^{*} (\fg , K;
C^{\infty} (\Gamma \backslash G)\otimes E).
\end{equation}
We can use \eqref{eqn:gKeqn} to identify important subspaces of the
cohomology.

For instance consider the space $L^{2} (\gbg)$.  Let $L^{2}_{\disc}
(\gbg )\subset L^{2} (\gbg)$ be the discrete part, and let
$L^{2}_{\cusp } (\gbg)$ be the subspace of cusp forms
\cite{harish-chandra}.  We have inclusions
\[
L^{2}_{\cusp} (\gbg)^{\infty } \hookrightarrow L^{2}_{\disc}
(\gbg)^{\infty} \hookrightarrow C^{\infty} (\gbg), 
\]
where the $\infty$ indicates that we take the subspaces of smooth
vectors.  These induce a map
\[
H^{*} (\fg , K; L^{2}_{\cusp} (\gbg)^{\infty }\otimes E)\rightarrow
H^{*} (\fg , K; C^{\infty} (\gbg)\otimes E)
\]
that is in fact injective.  The image $H_{\cusp}^{*} (\Gamma ; E)
\subset H^{*} (\Gamma ; E)$ is called the \emph{cuspidal cohomology}.
In a certain sense, this is the most interesting constituent 
of the cohomology.

Thus $H^{*} (\Gamma ; E)$ contains a subspace corresponding to certain
cuspidal automorphic forms.  It is natural to ask if the rest of the
cohomology can be built from automorphic forms as well.  More
precisely, let
\[
A (\Gamma , G) \subset C^{\infty} (\gbg)
\]
be the subspace of automorphic forms \cite{harish-chandra}, that is,
$A (\Gamma , G)$ is consists of functions of moderate growth that are
right $K$-finite and left $Z (\fg)$-finite, where $Z (\fg)$ denotes
the center of the universal enveloping algebra $U(\fg)$.  Then we have
the following theorem of Franke \cite{franke}, which verifies a
conjecture of Borel:

\begin{theorem}
\thmcite{franke}
The inclusion $A (\Gamma , G)\rightarrow C^{\infty} (\gbg)$ induces an
isomorphism
\[
H^{*} (\fg , K; A (\Gamma , G)\otimes E)\rightarrow  H^{*} (\fg , K; C^{\infty} (\gbg)\otimes E).
\]
\end{theorem}

Thus we can think of $H^{*} (\Gamma ; E)$ as being a concrete
realization of certain automorphic forms, namely those whose
associated automorphic representations have nonvanishing $(\fg ,
K)$-cohomology.  The unitary representations with this property were
classified by Vogan--Zuckerman \cite{vog.zuck}.  In general the
automorphic forms contributing to the cohomology are far from typical,
in the sense that most automorphic forms cannot appear in the
cohomology.  Nevertheless, those that do are very interesting, and are
expected to have direct connections with arithmetic
(cf.~\S\ref{s:arithmetic}).

\subsection{}\label{ss:modularforms} A prototypical example of the
connection between cohomology and automorphic forms is that of $\G =
\SL_{2}$.  In this case we have $G=\SL_{2} (\R)$, $K=\SO (2)$, and the
symmetric space $X$ is the upper halfplane.  Let $\Gamma = \Gamma_{0}
(N) \subset \SL_{2} (\Z)$, the subgroup of matrices upper triangular
modulo $N$.  Let $E_{k}$ be the $k$-dimensional complex representation
of $G$, say on the vector space of degree $k-1$ homogeneous complex
polynomials in two variables.  By work of Eichler and Shimura, we have
\[
H^{1} (\Gamma ; E_k) \simeq S_{k+1} (\Gamma )\oplus
\overline{S}_{k+1} (\Gamma) \oplus \Eis_{k+1} (\Gamma ),
\]
where $S_{k+1}$ is the space of holomorphic weight $k+1$ modular
forms, $\Eis_{k+1}$ is the space of weight $k+1$ Eisenstein series,
and the bar denotes complex conjugation.  Thus holomorphic modular
forms of weights $\geq 2$ appear in the cohomology; Maass forms and
weight $1$ holomorphic forms do not.

\def\vcd{\nu}

\subsection{}
In general the quotient $\Gamma \backslash X$ has cohomology in many
degrees. Which cohomology groups are the most interesting?

Clearly $H^{i} (\Gamma \backslash X; \E)=0$ if $i>\dim \Gamma
\backslash X = \dim X$.  But one actually knows that sometimes the
cohomology vanishes in degrees less than this.  More precisely, let $q
= q (\G)$ be the $\Q$-rank of $\G$, which by definition is the
dimension of a maximal torus split over $\Q$.  Then we have the
following theorem of Borel--Serre:

\begin{theorem}
\thmcite{borel.serre}
For all $\Gamma$ and $E$ as above, 
we have $H^{i} (\Gamma ; E) = 0$ if $i>\dim X - q$.
\end{theorem}

The number $\vcd = \dim X - q$ is called the \emph{virtual
cohomological dimension}.

One also knows that the cuspidal cohomology does not necessarily
appear in every cohomological degree.  In fact, for $\SL_{n}$ one can
show that $H^{i}_{\cusp} (\Gamma ; E) =0$ unless the degree $i$ lies
in a small interval about $(\dim X)/2$ \cite{li.schwermer.survey}.  As
we shall see, this makes computations much more difficult to perform
for large $n$.  Table \ref{dimtable} shows the cuspidal range for
subgroups of $\SL_{n} (\Z)$ for $n\leq 9$ \cite{schwermer.sln}.

\begin{table}[htb]
\begin{center}
\begin{tabular}{|c||c|c|c|c|c|c|c|c|}
\hline
$n$&2&3&4&5&6&7&8&9\cr
\hline
\hline
$\dim X$&2 & 5 & 9 & 14 & 20 & 27 & 35 & 44 \cr
$\vcd (\Gamma) $ & 1 & 3 & 6 & 10 & 15 & 21 & 28 & 36\cr
top degree of $H^{*}_{\text{cusp}} $& 1 & 3 & 5 & 8 & 11 & 15 & 19& 24 \cr
bottom degree of $H^{*}_{\text{cusp}} $& 1 & 2 & 4 & 6 & 9 & 12 & 16 & 20 \cr
\hline
\end{tabular}
\end{center}
\medskip
\caption{The virtual cohomological dimension and the cuspidal range
for subgroups of $\SL_{n} (\Z)$\label{dimtable}.}
\end{table}

\section{Connections with arithmetic geometry}\label{s:arithmetic}

The groups $H^{*} (\Gamma ; E)$ have an action of the \emph{Hecke
operators}, which are endomorphisms of the cohomology associated to
certain finite index subgroups of $\Gamma$.  The eigenclasses of these
operators reveal the arithmetic information in the cohomology.  We
review this now.

\subsection{}We begin with Galois representations.  Let $\G =
\SL_{n}/\Q$ and let $\Gamma = \Gamma_{0} (N)$ be as in the
introduction.  Let $\Gal (\bar \Q / \Q)$ be the absolute Galois group
of $\Q$.  Let $\rho \colon \Gal(\bar \Q / \Q)\rightarrow GL_{n}
(\Q_{p})$ be a continuous semisimple Galois representation unramified
outside $pN$.  For any prime $l$ not dividing $pN$ let $\Frob_{l}$ be
the Frobenius conjugacy class over $l$.  Then we can consider the
characteristic polynomial
\[
\det (1-\rho (\Frob_{l})T) \in \Q_{p}[T].
\]

\subsection{}
On the cohomology side, for each prime $l$ not dividing $N$ we have
Hecke operators $T (l,k)$, $k=1,\dotsc ,n-1$.  These operators
generalize the classical operator $T_{l}$ on modular forms; for an
exposition of these operators and the structure of the algebra they
generate, see \cite[Ch.~3]{shimura.book}. If $\xi$
is a simultaneous eigenclass for these operators, define the
\emph{Hecke polynomial}
\[
H (\xi) = \sum_{k} (-1)^{k}l^{k (k-1)/2} a (l,k) T^{k} \in \C [T].
\]
where $a (l,k)$ is the eigenvalue of $T (l,k)$.\footnote{Note that the
Hecke polynomial $H (\xi)$ is normalized such that if $H (\xi) $ were
used as a local factor for the $L$-function of the automorphic
representation $\pi_{\xi}$ attached to $\xi$, via the substitution
$T=l^{-s}$, then the functional equation of the $L$-function would
have the form $L(s,\pi_{\xi}) = L (n-s, \tilde{\pi}_{\xi})$, where the
tilde denotes contragredient.  This should be kept in mind when
studying the examples in \S\ref{s:results}.}

\subsection{} Now fix an isomorphism $\varphi \colon \C \rightarrow \bar \Q_{p}$.  Here is
one way to express the conjectural connection between cohomology and
arithmetic:

\begin{conjecture}\label{conj:motive}
For any Hecke eigenclass $\xi$ of level $N$, there is a Galois
representation $\rho \colon \Gal (\bar \Q /\Q)\rightarrow \GL_{n}
(\Q_{p})$ unramified outside $pN$ such that for every prime $l$ not
dividing $pN$, we have
\[
\varphi (H (\xi)) = \det (1-\rho (\Frob_{l})T).
\]
\end{conjecture}

This conjecture is a translation to the setting of the cohomology of
arithmetic groups of a more general conjecture of Clozel
\cite[Conjecture 4.5]{clozel.motives}, which precisely predicts the
cuspidal automorphic representations that should be attached to
motives.  This is the conjecture that we are ultimately testing.

We are primarily interested in Conjecture \ref{conj:motive} in the
case of eigenclasses that are not \emph{essentially selfdual}. (By
definition, a Hecke eigenclass is essentially selfdual if the
associated automorphic representation $\pi$ satisfies $\pi \simeq
\widetilde{\pi} \otimes \chi$, where the tilde denotes contragredient
and $\chi$ is a $1$-dimensional automorphic representation.)  This is
because in many cases one knows how to attach motives to essentially
selfdual eigenclasses by realizing the motive in the \'etale
cohomology of a Shimura variety (see for instance
\cite{clozel.galois}, which treats certain selfdual classes).  For
eigenclasses that are not essentially selfdual it is unknown in
general how to do this.  In our case, for $n\geq 3$ the symmetric
space $\Gamma \backslash X$ is not an algebraic variety, and it is
completely unclear how to directly connect the Galois group with Hecke
eigenclasses.

\section{Techniques}\label{s:techniques}

In this section, we describe the topological techniques we use to
compute the cohomology and the action of the Hecke operators.

\subsection{}\label{ss:4.1} For modular forms, i.e.~the cohomology of
subgroups of $\SL_{2} (\Z)$, we can use \emph{modular symbols} to
perform computations.  Recall that the symmetric space $X$ in this
case is the upper halfplane.  Let $X^{*} = X\cup \Q \cup \{i\infty
\}$, the standard partial compactification of $X$ obtained by
adjoining cusps; equip $X^{*}$ with the Satake topology.  Given two cusps
$q_{1}, q_{2}\in X^{*}\smallsetminus X$, we can form the ideal
geodesic in $X$ from $q_{1}$ to $q_{2}$ and can look at its image in
$\Gamma \backslash X^{*}$.  This gives a relative homology class
\begin{equation}\label{eqn:mstocohomology}
[q_{1}, q_{2}]\in H_{1} (\Gamma \backslash X^{*}, \text{cusps}; \C)
\simeq H^{1} (\Gamma \backslash X; \C).
\end{equation}

Let $V$ be the complex vector space generated by the symbols $[q_{1},
q_{2}]$.  Let $S$ be $V$ modulo the 
subspace spanned by all elements of form $[q_{1}, q_{2}] + [q_{2},q_{1}]$.
The group $\Gamma$ acts on $V$ by
\[
\gamma \cdot [q_{1}, q_{2}] = [\gamma q_{1}, \gamma q_{2}].
\]
This action descends to $S$, and we let $S_{\Gamma}$ be the quotient
of coinvariants.  The space $S_{\Gamma}$ is called the space of
\emph{modular symbols}.  We obtain a map
\begin{equation}\label{eqn:msiso}
S_\Gamma \rightarrow H^{1} (\Gamma ; \C)
\end{equation}
induced by \eqref{eqn:mstocohomology} whose kernel is easily
determined: it is spanned by sums of the form
\begin{equation}\label{eqn:threeterm}
[q_{1},q_{2}] + [q_{2}, q_{3}] + [q_{3}, q_{1}].
\end{equation}
Let $M_{\Gamma}$ be the quotient of $S_{\Gamma}$ by the subspace
generated by \eqref{eqn:threeterm}; this space is isomorphic to $H^{1}
(\Gamma ; \C)$.  

Thus modular symbols provide a combinatorial model for $H^{1} (\Gamma
; \C)$.  One can also enrich modular symbols by adding extra data to
account for the representation $E_{k}$ from \S\ref{ss:modularforms} to
compute with modular forms of higher weights; see \cite{stein.book}
for more details.

\subsection{} The Hecke action on the cohomology can be encoded in an
action on the modular symbols.  Suppose we are given a Hecke operator
$T_{l}$.  The construction of the Hecke algebra guarantees that we can
find finitely many matrices $\{\gamma_{i}\}\subset \GL_{2} (\Q) \cap
M_{2} (\Z)$ such that, if $\xi = [q_{1}, q_{2}] \in S_{\Gamma}$, then
\begin{equation}\label{eqn:tl}
T_{l}\xi = \sum_{i} [\gamma_{i}q_{1}, \gamma_{i}q_{2}].  
\end{equation}
It is easy to verify that this is an action of the Hecke operators
on $S_{\Gamma}$ that is compatible with the isomorphism \eqref{eqn:msiso}.

\subsection{} Hence we have a concrete model for the cohomology space
$H^{1} (\Gamma ; \C)$, as well as a Hecke action on our model.  But
unfortunately this is not good enough for practical computation of
Hecke eigenvalues.

The first problem is that the space $S_{\Gamma }$ is
infinite-dimensional.  Since $H^{1} (\Gamma ; \C)$ is
finite-dimensional, one might expect to be able to identify a
finite-dimensional subspace of $S_{\Gamma}$ that maps onto
$M_{\Gamma}$.  A very convenient subspace with this property is
provided by the \emph{unimodular symbols}.  These are the images in
$S_{\Gamma}$ of the $\SL_{2} (\Z)$-translates of $[0,i\infty]$.  It is
easy to see that these images generate a finite-dimensional subspace
that is still large enough to surject onto the cohomology.  For
instance, if $\Gamma = \Gamma (N)$, the principal congruence subgroup
of matrices congruent to the identity mod $N$, then on $X (N) = \Gamma (N)
\backslash X^{*}$ the ideal geodesics giving the unimodular symbols
become the edges of a highly-symmetric triangulation.  The unimodular
subspace of $S_{\Gamma}$ is just the first chain group attached to
this triangulation, and the relations \eqref{eqn:threeterm} correspond
to the boundaries of triangles.

\subsection{} But now one encounters another problem.  The Hecke operators
act on modular symbols, and the unimodular symbols span the
cohomology, but the Hecke operators do not preserve the subspace of
unimodular symbols.  In fact, no finite set of modular symbols modulo
$\Gamma$ admits a Hecke action.

This is easy to see; one can attach a integer ``determinant'' $n
(\xi)$ to a modular symbol $\xi = [q_{1},q_{2}]$ as follows.  First
write $q_{i}=a_{i}/b_{i}$, $i=1,2$, where the rational numbers are in
reduced terms.  Then define $n (\xi) = |a_{1}b_{2}-a_{2}b_{1}|$.  The
integer $n (\xi)$ is well-defined modulo the defining relations for
$S_{\Gamma}$.  We have $n (\xi)= 1$ if and only if $\xi $ is
unimodular.\footnote{This is of course why they are called
{unimodular}.}  Now if one applies $T_{l}$ to a modular symbol
$\xi$ with $n (\xi) = d$, then in general the modular symbols on the
right of \eqref{eqn:tl} will have determinant $ld$.  Since we can take
$l$ to arbitrarily large, no finite spanning set can be preserved by
all the $T_{l}$.

Nevertheless one can circumvent this difficulty.  There is an
algorithm---due to Manin \cite{manin}---that enables one to write any
modular symbol as a linear combination of unimodular symbols.
Essentially Manin's algorithm is nothing more than the continued
fraction algorithm.  Conjugating by $\SL_{2} (\Q)$ we may assume the
nonunimodular symbol $\xi$ has the form $[0,q]$.  Then $\xi$ is taken
to a linear combination of the symbols $\xi_{i} = [a_{i}/b_{i},
a_{i+1}/b_{i+1}]$, where the $a_{i}/b_{i}$ are the convergents in the
continued fraction expansion of $q$.  Standard facts from elementary
number theory imply that the $\xi_{i}$ are unimodular.  This enables
us to compute the effect of a Hecke operator on a basis of unimodular
symbols, which allows us compute the Hecke eigenvalues and
eigenclasses.

\subsection{} Now we consider $n>2$.  First of all, there is an
analogue of the modular and unimodular symbols, and they provide a
model for $H^{\vcd} (\Gamma ; \C)$.  The space of modular symbols is
generated by $n$-tuples of cusps $[q_{1},\dotsc ,q_{n}]$ modulo
relations, where now cusp means a minimal boundary component in a
certain Satake compactification $X^{*}$ of $X$ \cite{satake2}.  There
is an analogue of Manin's algorithm, due to Ash--Rudolph
\cite{ash.rudolph}, that enables us to compute the Hecke action on
$H^{\vcd} (\Gamma ; \C)$.

The problem now is that \emph{usually} $H^{\vcd}_{\cusp} =0$.  In
fact, looking at Table \ref{dimtable}, we see that $\vcd$ lies in the
cuspidal range only for $\SL_{2}$, $\SL_3$.  Thus in these cases we
can use modular symbols, together with the Ash--Rudolph algorithm, to
compute Hecke eigenclasses.  For $\SL_{3}$ this was done originally by
Ash--Grayson--Green \cite{ash.grayson.green}, and later by
van~Geemen--Top~\cite{vgt1, four.dutch}, as well as some other groups.
As a result of these investigations, we have lots of convincing
evidence that cuspidal cohomology is related to arithmetic geometry
\cite{apt, vgt1, exp.ind}.

\subsection{} Our case of interest is $n=4$.  Here there are several
tools that allow us to overcome the problems above.  Most of these
tools work for all $n$: in principle we can compute the 
cohomology for all $n$, even though this computation quickly gets
impractical as $n$ increases.  On the other hand, we do not know if
the Hecke operator algorithm will generalize to allow computation of
the Hecke eigenvalues in the cuspidal range for $n>4$.

We can construct an explicit homological complex $S_{k}$, $k\geq 0$,
the \emph{sharbly complex}, that computes $H^{*} (\Gamma ; \C)$.  This
complex, defined by Ash \cite{ash.sharb}, generalizes the modular
symbols, but now we have chain groups in all nonnegative degrees.  The
sharbly complex has a $\Gamma$-action, and by the Borel--Serre duality
theorem \cite{borel.serre} the homology of the complex of coinvariants
$(S_{*})_{\Gamma}$ gives the cohomology of $\Gamma$.  For example, in
the case of $\SL_{2}$, the space $S_{0}$ coincides with the space $S$
from \S\ref{ss:4.1}, and the image of the boundary $S_{1}\rightarrow
S_{0}$ consists of sums of the form \eqref{eqn:threeterm}.  Passing to
the coinvariants we obtain our presentation of $H^{1} (\Gamma ; \C)$

The sharbly complex admits a Hecke action, and as before there are
infinitely many sharblies modulo $\Gamma$.  The analogue of the
unimodular symbols is a certain subcomplex of \emph{reduced sharblies}
that is finite mod $\Gamma$.  This subcomplex can be constructed using
the well-rounded retract of Ash \cite{ash.wrr}, or equivalently
Voronoi's work on perfect quadratic forms \cite{voronoi1}.

To compute the Hecke operators, we need an analogue of the Manin trick
to compute with $H^{\vcd -1} (\Gamma ; \C)$.  This is actually much
more subtle since we work below the cohomological dimension.  The
technique we have developed to do this has been described in several
places \cite{experimental, gy1, stein.book}, so we only give the idea.
The main point is that to compute the cohomology $H^{\vcd -1} (\Gamma
; \C)$ one works with sharbly cycles consisting of $n$-simplices glued
together along ``submodular symbols.''  (In the case of $\SL_{2}
(\Z)$, such cycles are built of geodesic triangles with usual modular
symbols as faces.)  To compute the Hecke action, one must move a
general sharbly cycle to a sum of reduced cycles by simultaneously
applying the Ash--Rudolph algorithm to each of its submodular symbols.
For more details we refer to the above references.

\begin{remark}
The technique described in \cite{experimental} to compute the Hecke
operators on the cohomology of subgroups of $\SL_{4} (\Z)$ is also useful in other
contexts.  One needs a linear group where the cuspidal range goes up
to one below the cohomological dimension.  We have explored this in
some cases:
\begin{itemize}
\item $\G = R_{F/\Q} (\GL_{2})$, where $F$ is real quadratic
\cite{gy1}.  This is the \emph{Hilbert modular case}; the cohomology
classes are related to weight $(2,2)$ Hilbert modular forms over $F$.
Note that we use $\GL_{2}$ instead of $\SL_{2}$.  We do this because
the $\GL_{2}$-symmetric space $X_{\GL_{2}}$ can be represented as a
certain cone of positive-definite binary quadratic forms modulo
homotheties.  One can then apply the generalization of \Vor's theory
of perfect forms, see for example
\cite{koecher1, koecher2, amrt, ash.oldpaper} to get a cellular
decomposition of $X_{\GL_{2}}$.  However this does not affect the
results of the computation in a serious way, because the associated
locally symmetric space $\Gamma \backslash X_{\GL_{2}}$ is a circle
bundle over a Hilbert modular surface.
\item $\G = R_{F/\Q} (\GL_{2})$, where $F$ is complex quartic
\cite{ghy}.  Again we work with the $\GL_{2}$-symmetric space.
\end{itemize}
\end{remark}

\section{Results}\label{s:results}

We have computed $H^{5} (\Gamma_{0} (N); \C)$ for $N$ prime and $\leq
211$, and for composite $N$ up to $52$.  Our biggest computation
involved computing the Smith normal form of a matrix of size $ 845712
\times 3277686$ ($N=211$).  Such computations cannot be done by merely
using existing numerical linear algebra software, because of the
special needs of our computation.  For instance, instead of working
over $\C$ we actually work over a large finite field.  We also need to
compute change of basis matrices to be able to compute Hecke
operators.  The specific numerical techniques we developed are
described in detail in \cite[\S2]{agm3}.

Our computation results can be summarized as follows:
\begin{itemize}
\item In the range of our computations, we found no nonselfdual
cuspidal classes.  We know of no reason why they should not exist for
larger~$N$, but no one has proven their existence.\footnote{At the
conference, Neil Dummigan explained how some of his recent work
suggests that there should exist nonselfdual classes on $\SL_{4}$ at
level $1$ and high weight. }

\item We found Eisenstein classes (\S\ref{ss:eisenstein}) attached to
weight 2 and weight 4 modular forms (\S\ref{ss:twoandfour}).
\item We found Eisenstein classes attached to $\SL_{3}$ cuspidal
cohomology (\S\ref{ss:sl3}).  
\item We found selfdual cuspidal classes that are apparently functorial
lifts of Siegel modular forms (\S\ref{ss:siegel}).
\end{itemize}

For $N$ prime we believe that our description of the Eisenstein
cohomology is complete.  However, if $N$ is not prime then there are
other Eisenstein classes, as mentioned already in \cite{agm1}.

\subsection{Eisenstein cohomology}\label{ss:eisenstein} 
We return for the moment to a general setting.  Our goal is to explain
the notion of Eisenstein cohomology, due to Harder \cite{harder.icm}.
We use the notation of \S\ref{ss:notation}, so that $X$ is a global
symmetric space and $\Gamma$ is an arithmetic group.

Let $\bar X$ be the partial compactification of $X$ due to
Borel--Serre \cite{borel.serre}.  This is an enlargment of $X$
that adds additional spaces at infinity.  The quotient $\Gamma
\backslash \bar X$ is called the \emph{Borel--Serre compactification}
of $\Gamma \backslash X$.  If $\Gamma$ is torsion-free, then $\Gamma
\backslash \bar X$ has the structure of a manifold with corners.
Otherwise, $\Gamma \backslash \bar X$ is an orbifold with corners.
Let $\partial (\Gamma \backslash \bar X) = \Gamma \backslash \bar
X\smallsetminus \Gamma \backslash X$ be the boundary.

The homotopy equivalence $\Gamma \backslash \bar X \rightarrow \Gamma
\backslash X$ induces an isomorphism
\[
H^{*} (\Gamma \backslash \bar X; \C) \simeq H^{*} (\Gamma \backslash X; \C),
\]
and the inclusion $\partial (\Gamma \backslash \bar X) \hookrightarrow
\Gamma \backslash \bar X$ induces a restriction map
\[
H^{*} (\Gamma \backslash \bar X; \C) \rightarrow H^{*} (\partial
(\Gamma \backslash \bar X); \C ).
\]
The \emph{Eisenstein classes} are certain classes that form a basis of
the image of a splitting of this restriction map.

\subsection{Weights $2$ and $4$}\label{ss:twoandfour} Now let $\Gamma
= \Gamma_{0} (N)\subset \SL_{4} (\Z)$, where $N$ is prime.  In $H^{5}
(\Gamma ; \C)$ we found Eisenstein classes corresponding to weight $2$
and weight $4$ holomorphic modular forms of level $N$.

Let $f$ be a weight $2$ newform of level $N$.  Then $f$ contributes to
$H^{5} (\Gamma ; \C)$ in two different ways, with the Hecke
polynomials
\[
(1-l^{2}T) (1-l^{3}T) (1-\alpha T + lT^{2})
\]
and 
\[
(1-T) (1-l T)(1-l^{2}\alpha T + l^{5}T^{2}),
\]
where $T_{l}f = \alpha f$.  
 
Now let $g$ be a weight $4$ newform.  Then whether or not $g$
contributes an Eisenstein class depends on the central special value
of the $L$-function of $g$.  Namely, if this special value vanishes,
then $g$ corresponds to an eigenclass with Hecke polynomial
\[
(1-lT) (1-l^{2}T) (1-\beta T + l^{3}T^{2}),
\]
where $T_{l}g=\beta g$.  This dependence on special values is typical
for Eisenstein cohomology.  We remark also that this particular
Eisenstein class is apparently a ``ghost class'' in the sense of
\cite{harder.1447}. That is, this class in $H^5(\Gamma ; \C)$ restricts
nontrivially to the boundary of the Borel--Serre compactification, but
it restricts trivially on each face of the boundary.

\subsection{$\SL_{3}$ cuspidal classes}\label{ss:sl3}

We also found cohomology classes corresponding to cuspidal cohomology
classes on subgroups of $\SL_{3} (\Z)$.  
These cohomology classes were originally computed by
Ash--Grayson--Green \cite{ash.grayson.green}.

Let $\eta \in H_{\cusp}^{3} (\Gamma_{0}^{*} (N); \C)$ be a cuspidal
cohomology class on $\Gamma^{*}_{0} (N) \subset \SL_{3} (\Z)$.  Let
$\gamma , \gamma '$ be the Hecke eigenvalues of this class with
respect to the Hecke operators $T (l,1), T (l,2)$ (computed of course
when $n=3$.).  Then $\eta$ contributes in two different ways, with the
Hecke polynomials
\[
(1-l^{3}T) (1-\gamma T + l \gamma 'T^{2} - l^{3}T^{3})
\]
and 
\[
(1-T) (1-l\gamma T + l^{3} \gamma 'T^{2} - l^{6}T^{3}).
\]

\subsection{Siegel modular forms}\label{ss:siegel} Finally, we
describe the contributions to the cohomology coming from Siegel
modular forms.  For more background on Siegel modular forms we refer
to \cite{vdg, ibuk2}.

Let $K (p)$ be the \emph{paramodular group} of prime
level, which by defintion is the subgroup
\[
K (p) = \left(\begin{array}{cccc}
\Z & p\Z & \Z & \Z \\
\Z & \Z &\Z &p^{-1}\Z \\
\Z &p\Z &\Z&\Z \\
p\Z &p\Z&p\Z &\Z   
\end{array} \right) \subset \Sp_{4} (\Q).
\]
Let $S^{3} (p)$ be the space of weight three paramodular forms.  Note
that all such forms are cuspforms; there are no Eisenstein series.
This space contains the subspace $S_{\grit}^{3} (p)$ of
\emph{Gritsenko lifts}, which are lifts from certain weight $3$ Jacobi
forms to $S^3 (p)$ \cite{grit}.  Let $S^{3}_\nongrit (p)$ be the Hecke
complement to $S^{3}_{\grit}(p)$ in $S^{3} (p)$.  The forms in
$S^{3}_\nongrit (p)$ will be those that appear in the cohomology.

The space of cuspidal paramodular forms at prime level is known pretty
explicitly.  First we have a dimension formula due to Ibukiyama
\cite{ibuk1, ibuk2}.

Let $\kappa(a)$ be the Kronecker symbol $(\frac{a}{p})$.  Define
functions $f, g\colon \Z \rightarrow \Q$ by
\[
f (p) = \begin{cases}
2/5&\text{if $p\equiv 2, 3 \bmod 5$,}\\
1/5&\text{if $p=5$,}\\
0&\text{otherwise},
\end{cases}
\]
and 
\[
g (p) = \begin{cases}
1/6&\text{if $p\equiv 5\bmod 12$,}\\
0&\text{otherwise}.
\end{cases}
\]
\begin{theorem}[Ibukiyama]
For $p$ prime 
we have $\dim S^{3} (2) = \dim S^{3} (3) = 0$.  For
$p\geq 5$, we have 
\begin{align*}\label{eq:ibuk}
\dim S^{3} (p) &= (p^{2}-1)/2880 \\
&+ (p+1) (1-\kappa (-1))/64 
+ 5 (p-1) (1+\kappa (-1))/192 \\
&+ (p+1) (1-\kappa (-3))/72 + (p-1)
(1+\kappa (-3))/36 \\
&+ (1-\kappa (2))/8+f (p)+g (p) -1.
\end{align*}
\end{theorem}

Using this one can easily compute the dimension of $S_{\nongrit}^{3} (p)$.

Next, Poor and Yuen \cite{py1, py2} have developed a technique to compute Hecke
eigenvalues for weight three paramodular forms, and in particular for
the forms in $S_{\nongrit}^{3} (p)$.   Using Ibukiyama's formula and
data supplied to us by Poor and Yuen, we find the following:

\begin{itemize}
\item For all prime levels $N$, the dimension of the subspace of $H^{5} (\Gamma_{0}
(N);\C) $ not accounted for by the Eisenstein classes above matches $2
\dim S^{3}_{\nongrit} (N)$.
\item In cases where we have computed the Hecke action on this
subspace, we find full agreement with the data produced by
Poor--Yuen.  That is, the Hecke polynomials of the paramodular forms
exactly match the polynomials we compute for an eigenbasis of the
complement to the Eisenstein classes.
\end{itemize}

\subsection{Future problems}
We conclude by indicating some open questions and prospects for future work.

\begin{itemize}
\item Can one prove that the Eisenstein classes we found actually
occur in the cohomology?  In principle, one should be able to apply
standard techniques from Eisenstein cohomology to answer this
question, although it appears at the moment that current knowledge is
insufficient to provide a proof.
\item Can one prove the lift from paramodular forms to the cohomology?
We have not tried to investigate this question at all.  Perhaps it is
not difficult to solve using known cases of functoriality in the literature.
\item Study the Hecke action on torsion classes in the cohomology.
This is currently under investigation in joint work with Ash and McConnell.
\end{itemize}
\bibliographystyle{amsplain_initials}
\bibliography{rims}

 \end{document}